\documentclass[twocolumn,a4paper,10pt]{article}
\usepackage{amssymb,amsmath}
\usepackage[applemac]{inputenc}
\usepackage{xypic}
\xyoption{all}
\usepackage{graphicx}

\newtheorem{theorem}{Theorem.}

\newcommand\V{\bigvee}
\newcommand\Max{\operatorname{Max}}
\newcommand\ie{i.e.}
\newcommand\eg{e.g.}
\newcommand\st{\mid}
\newcommand\opens{\operatorname{\mathcal{O}}}
\newcommand\spectrum{\operatorname{\Sigma}}
\newcommand\topology{\operatorname{\Omega}}
\newcommand\groupoid{\operatorname{\mathcal{G}}}
\newcommand\downsegment{{\downarrow}}
\newcommand\isgrp{\Gamma}
\newcommand\lcc{\operatorname{{\mathcal L}^{\vee}}}
\newcommand\Germs{\operatorname{Germs}}
\newcommand\sheaf{\mathcal D}

\begin{document}

\title{Quantales as geometric objects: symmetry beyond groupoids?}
\author{Pedro Resende\thanks{Supported in part by FEDER
and FCT through CAMGSD.}
\vspace*{2mm}\\ \small\it Departamento de Matem{\'a}tica, Instituto
Superior T{\'e}cnico, \vspace*{-1mm}\\ \small\it Av. Rovisco Pais 1,
1049-001 Lisboa, Portugal}

\date{~}

\maketitle

\section{Introduction}

Modern mathematics has become pervaded by the idea that in order to cater for certain notions of \emph{symmetry}, in particular of a local nature, one needs to go beyond group theory, replacing groups by \emph{groupoids}. A nice survey of some implications of this idea in algebra, geometry, and analysis (as of 1996) can be found in \cite{Weinstein}. The same idea can be expressed in terms of \emph{pseudogroups}, which provide another generalization of groups, or, even more generally, by abstract \emph{inverse semigroups}. See the first chapters of \cite{Lawson} for motivations and a good historical account.

How do these two generalizations relate to each other? Many differences and similarities are illustrated by various constructions back and forth between groupoids and inverse semigroups \cite{Lawson,Paterson}, and by their relations to operator algebras \cite{Connes,Paterson,Renault}. One recurrent aspect is that inverse semigroups are closely related to a particular notion of \emph{topological groupoid}, namely to \emph{étale groupoids}. A fundamental reason for this can be singled out in the form of a precise correspondence between these two concepts, bearing a close resemblance to the equivalence between local homeomorphisms and sheaves: to each inverse semigroup $S$ of a certain kind (the analogue of a sheaf) we associate the groupoid of \emph{germs} of $S$, which is an étale groupoid $\Germs(S)$ (the analogue of a local homeomorphism); and to each étale groupoid $G$ we associate an inverse semigroup of ``sections'' $\isgrp(G)$, such that $S\cong\isgrp(\Germs(S))$ and $G\cong\Germs(\isgrp(G))$. We provide a brief description of this correspondence in \S\ref{sec:gi}.

In order to extend the correspondence beyond étale groupoids we need a more general type of semigroup. One good candidate is the notion of \emph{quantale} \cite{Mul:enc,Rosenthal}, of which the most general definition is that of a \emph{sup-lattice ordered semigroup} (\S\ref{sec:q}). For instance, the topology of any étale groupoid is closed under pointwise multiplication of open sets, hence being a quantale. There is an algebraic characterization of the quantales obtained in this way, through which a correspondence between étale groupoids and quantales is obtained \cite{egatq}, matching that of inverse semigroups and étale groupoids (\S\ref{sec:qg}). But we can also go beyond étale groupoids because slightly more general classes of quantales provide characterizations of more general groupoids, such as open groupoids (\S\ref{sec:og}).

The purpose of this paper, which is to be regarded partly as a research announcement, is to highlight some aspects of the interplay between quantales, inverse semigroups, and groupoids. Many of the results mentioned have not yet been presented (some are joint work) and will appear in detail elsewhere.

\section{Groupoids}\label{sec:d}

A \emph{groupoid} $G$ is a small category in which every morphism is invertible, or, equivalently, a pair of sets $G_0$ (of \emph{units}) and $G_1$ (of \emph{arrows}) equipped with structure maps satisfying appropriate axioms, 
\[
\xymatrix{
G_1\times_{G_0}G_1\ar[r]^-m&G_1\ar@(ru,lu)[]_i\ar@<1.2ex>[rr]^r\ar@<-1.2ex>[rr]_d&&G_0\ar[ll]|u
}
\]
where the multiplication map $m$ is defined on the set $G_1\times_{G_0}G_1$ of \emph{composable pairs of arrows}:
\[G_1\times_{G_0}G_1=\{(x,y)\in G_1\times G_1\st r(x)=d(y)\}\;.\]
The above diagram makes sense in any category with pullbacks, for then the ``object of composable pairs of arrows'' can be defined, and in particular it makes sense in the category of topological spaces, where it gives us topological groupoids. A topological groupoid is \emph{open} if $d$ (equivalently, all the structure maps) is open, and \emph{étale} if $d$ (equivalently, all the structure maps) is a local homeomorphism. A related notion is that of \emph{r-discrete} groupoid \cite{Renault}, which in the applications considered in \cite{Renault} (due to the presence of a suitable measure) is the same as an étale groupoid.

In the category of \emph{locales} (see \S\ref{sec:l} below) there are pullbacks and thus we can define groupoids in it, obtaining the notion of \emph{localic groupoid}. Many topological definitions can be easily transferred to locales, in particular open maps and local homeomorphisms. \emph{Open} localic groupoids and \emph{\'{e}tale} localic groupoids are defined accordingly. From any localic groupoid $G$ there is a canonically associated topological groupoid $\spectrum (G)$ (its \emph{spectrum}), whose spaces of arrows and units are respectively the spectra of the locales $G_1$ and $G_0$.

\section{Inverse semigroups}

An \emph{inverse semigroup} (see \cite{Lawson}) is a semigroup $S$ equipped with an involution $s\mapsto s^{-1}$ that satisfies $ss^{-1}s=s$ and such that all the idempotents commute. The set of idempotents $E(S)$ is a semilattice. As an example, if $X$ is a topological space with topology $\topology(X)$ there is an inverse semigroup $\Gamma(X)$ (an example of a pseudogroup) whose elements are the homeomorphisms $h:U\to V$ with $U,V\in\topology(X)$. The idempotents are the identities on open sets, and thus $E(S)\cong\topology(X)$. By the Vagner-Preston representation theorem, every inverse semigroup is, up to isomorphism, contained in a pseudogroup.

Every inverse semigroup $S$ has a \emph{natural order} defined by
$s\le t$ if $s=ft$ for some idempotent $f$, which in $\Gamma(X)$ coincides with the restriction order. Two elements $s,t\in S$ are \emph{compatible} if $st^{-1}, s^{-1}t\in E(S)$, and $S$ is \emph{complete} if every set  $X\subseteq S$ whose elements are pairwise compatible has a supremum, or join, $\V X$ in $S$. The representation given by the Vagner-Preston theorem does not necessarily preserve joins. This follows immediately from the fact that every locale (see \S\ref{sec:l} below) is an inverse semigroup (all the elements are idempotents), but not every locale is spatial.

\section{Étale groupoids and inverse semigroups}\label{sec:gi}

Let $G$ be a topological étale groupoid. For simplification let us assume that $G_0$ is a sober topological space (for instance a Hausdorff space, as is often assumed in applications --- see \S\ref{sec:l}). A continuous local section $s:U\to G_1$ of $d$ is a \emph{local bisection} if $r\circ s$ is a homeomorphism onto its image. It is easy to see that the local bisections form a sheaf of sets on $G_0$. In addition, the set of all the local bisections has the structure of a complete inverse semigroup $S=\isgrp(G)$ for which $E(S)\cong\topology(G_0)$, where the multiplication is defined in terms of the multiplication of $G$ in a straightforward manner.

There is a converse to this construction. For each complete inverse semigroup $S$ whose semilattice of idempotents $E(S)$ is (isomorphic to) the topology $\topology(X)$ of a sober space $X$ we can define a sheaf of sets $\sheaf$ where for each ``open'' $U\in E(S)$ the set $\sheaf(U)$ contains the elements $s\in S$ such that $ss^{-1}=U$. If $V\le U$ are idempotents, the restriction map $\sheaf(U)\to\sheaf(V)$ is given by multiplication: $s\mapsto Vs$. The fact that this is a sheaf rather than just a presheaf is precisely equivalent to the completeness of $S$. Now the standard construction of a local homeomorphism from a sheaf (see, \eg, \cite[Ch.\ II.5]{MM}) gives us a space $\Lambda_\sheaf$ of ``germs'', along with a local homeomorphism $d:\Lambda_\sheaf\to X$. This is the domain map of an étale groupoid $\Germs(S)$ whose other structure maps are obtained from the inverse semigroup structure of $S$.

If $S=\isgrp(G)$ for an étale groupoid $G$, every local section of $d$ is ``locally a local bisection'', and thus $\Lambda_\sheaf$ is homeomorphic to the space of germs of local sections of $d:G_1\to G_0$, which is homeomorphic to $G_1$. It follows that both $G_1$ and the domain map $d:G_1\to G_0$ are recovered from $S$. It can be verified that the remaining structure maps of the groupoid $\Germs(S)$ agree with those of $G$, giving us an isomorphism $\Germs(\isgrp(G))\cong G$. Also, we have $\isgrp(\Germs(S))\cong S$.

Hence, a topological étale groupoid (with a sober space of units) is essentially ``the same'' as a complete inverse semigroup with a spatial locale of idempotents. These results can be generalized to unit spaces that are not sober.

\section{Quantales}\label{sec:q}

A \emph{unital involutive quantale} is a \emph{sup-lattice} (\ie, a partially ordered set in which every subset $X$ has a join $\V X$ --- and therefore also an infimum\footnote{Sup-lattices are complete lattices. The name ``sup-lattice'' is motivated \cite{JT} by thinking of joins as the first class operations, with meets being just derived. Accordingly, the homomorphisms of the category of sup-lattices are required to preserve only joins.}, or meet, $\bigwedge X$) equipped with an additional structure of involutive monoid (the involution is usually denoted by $a\mapsto a^*$, and the multiplicative unit by $e$), where the involution preserves joins, and so does the multiplication in each variable:
\begin{eqnarray*}
a({\V_i b_i}) &=& \V_i(a b_i)\\
({\V_i a_i})b&=&\V_i(a_i b)\\
({\V_i a_i})^*&=&\V a_i^*\;.
\end{eqnarray*}
A \emph{homomorphism}  $f:Q\to R$ of unital involutive quantales is a function that preserves joins, the multiplication, the multiplicative unit, and the involution:
\begin{eqnarray*}
f(\V S)&=&\V f(S)\\
f(a b)&=&f(a) f(b)\\
f(e_Q)&=& e_R\\
f(a^*)&=&f(a)^*\;.
\end{eqnarray*}
Quantales are ring-like structures, and there are corresponding notions of \emph{module}. A \emph{left module} over a unital quantale $Q$ is a sup-lattice $M$ equipped with a left action
$Q\times M\to M$ that preserves joins in each variable. For details on modules over involutive quantales see for instance \cite{Re04}.

\section{Locales}\label{sec:l}

To a large extent, the results with which we shall be concerned can be conveniently expressed in the language of \emph{locale theory}, of which we give a very basic outline (for details see \cite{StoneSpaces}).

By a \emph{locale} is meant a sup-lattice in which binary meets distribute over arbitrary joins:
\[x\wedge\V Y=\V_{y\in Y}(x\wedge y)\;.\]
(Hence, any locale is a unital involutive quantale with multiplication given by $\wedge$ and trivial involution $a^*=a$.)
The motivating example of a locale is the topology $\topology (X)$ of a topological space $X$, ordered by inclusion of open sets. By a \emph{map} $f:A\to B$ of locales is meant a \emph{homomorphism}
$f^*:B\to A$, \ie, a function that preserves arbitrary joins and finite meets (including the empty meet $1_B=\bigwedge\emptyset=\V B$):
\begin{eqnarray*}
f^*(\V S)&=&\V f^*(S)\\
f^*(a\wedge b)&=&f^*(a)\wedge f^*(b)\\
f^*(1_B)&=& 1_A\;.
\end{eqnarray*}
Again, the motivating example is the map of locales $\topology (X)\to\topology (Y)$  which is defined by the inverse image homomorphism
$f^{-1}:\topology (Y)\to\topology (X)$ of a continuous map $f:X\to Y$ of topological spaces.

A \emph{point} of a locale $A$ is defined to be a map of locales $p:\Omega\to A$ from the topology $\Omega$ of a discrete one point space to $A$. The \emph{spectrum} of a locale $A$ is the topological space $\spectrum (A)$ consisting of the points of $A$ with open sets of the form
\[U_a=\{p:\Omega\to A\st p^*(a)=1\}\;.\]
(This defines a functor $\spectrum$ from locales to topological spaces.)
The assignment $a\mapsto U_a$ is a surjective homomorphism of locales. $A$ is said to be \emph{spatial} if this is an isomorphism.

If a space $X$ is the spectrum of a locale then there is a homeomorphism $\spectrum(\topology(X))\cong X$. Spaces with this property are called \emph{sober} (\eg, any Hausdorff space). The category of sober spaces with continuous maps is equivalent to the category of spatial locales and their maps.

Locales are often important as replacements for the notion of topological space in a \emph{constructive} setting, by which is meant the ability to interpret definitions and theorems in an arbitrary topos. See \cite{J}. For instance, the locale $\operatorname{RIdl}(R)$ of radical ideals of a commutative ring $R$ can be regarded as the ``constructive Zariski spectrum'' of $R$ because $\spectrum(\operatorname{RIdl}(R))$ is (classically) homeomorphic to the usual space of prime ideals with the Zariski topology.

\section{Étale groupoids and quantales}\label{sec:qg}

Let $G$ be a topological étale groupoid. The fact that all the structure maps are local homeomorphisms implies two immediate facts: the unit space $G_0$ (rather, its image $u(G_0)$) is open in $G_1$; the pointwise product of any two open sets $U,V\in\topology(G_1)$ is an open set. This makes $\topology(G_1)$ a unital (and involutive) quantale. A topological groupoid is étale precisely if its topology has this property \cite[Th.\ 5.18]{egatq}.

The algebraic characterization of the unital involutive quantales that arise in this way has been done in \cite{egatq}, leading to a correspondence between localic étale groupoids and certain quantales. From a localic étale groupoid $G$ one obtains a quantale $\opens(G)$ through the localic analogue of the construction just described. For want of a better name, let us refer to such quantales as \emph{étale groupoid quantales}. Among other things they are also locales. The converse construction yields, from an étale groupoid quantale $Q$, a localic groupoid $\groupoid(Q)$ whose locale of arrows is $Q$, whose locale of units is $\downsegment e=\{a\in Q\st a\le e\}$, and such that \cite{egatq}
\begin{eqnarray}
G&\cong&\groupoid(\opens(G)) \label{eq1}\\
Q&\cong&\opens(\groupoid(Q))\;. \label{eq2}
\end{eqnarray}
The correspondence between inverse semigroups and topological étale groupoids can now be recast in terms of these quantales. Let $S$ be an \emph{abstract pseudogroup}, by which will be meant an inverse semigroup $S$ whose idempotents form a locale $E(S)$, and let us denote by $\lcc(S)$ the set of all the downwards closed sets of $S$ which are closed under the formation of all the existing joins of $S$.
\begin{theorem}
$\lcc(S)$ is an étale groupoid quantale. If in addition $E(S)$ is a spatial locale then $\lcc(S)$ is a spatial locale, and the (spectrum of) the groupoid of $\lcc(S)$ is the germ groupoid of $S$:
\[\spectrum(\groupoid(\lcc(S)))\cong\Germs (S)\;.\]
\end{theorem}
We are therefore provided with a generalization of the construction of germ groupoids to the localic setting.

The ``duality'' expressed by (\ref{eq1}) and (\ref{eq2}) does not extend to an equivalence of categories because the inverse image locale homomorphism $h^*$ of a map of localic groupoids $h:G\to H$ is not the same as a homomorphism of étale groupoid quantales $\opens(H)\to\opens(G)$ \cite{egatq}\footnote{A way around this would be to either expand or restrict the classes of morphisms under consideration. An analogous situation occurs with the equivalences of categories between inverse semigroups and inductive groupoids in \cite[Ch. 4, p.\ 114]{Lawson}.}. A consequence of this is that we are provided with an alternative category (a subcategory of quantales) whose objects are the étale groupoids. This category may be the right one to consider in some situations. For instance, the assignment $S\mapsto\lcc (S)$ is part of a left adjoint functor from abstract pseudogroups to étale groupoid quantales, and thus the identification of $\lcc (S)$ with a groupoid allows us to think of $\lcc (S)$ as being the ``universal'', or ``enveloping'', groupoid of $S$, with the proviso that the universality should be understood in the category of quantales rather than groupoids (in other words, roughly, it is the topology of the groupoid that is ``freely'' generated, rather than the groupoid itself).
Paterson's universal groupoid of an inverse semigroup $S$ \cite{Paterson} coincides with the groupoid of germs of a larger inverse semigroup $S'$, but the universality described in \cite[Prop.\ 4.3.5]{Paterson} is different.

In fact the adjunction just mentioned takes place between abstract pseudogroups and a larger category of quantales (the category of \emph{stable quantal frames} \cite{egatq}). The latter deserves attention because it has good properties, but besides the étale groupoid quantales it contains other quantales. These can be identified with involutive graphs that are almost étale groupoids, with the exception that their multiplication is ``fuzzy'' because it assigns to each composable pair of arrows an open set of arrows rather than just an arrow \cite[\S4.4]{egatq}. The usefulness of such a generalization in applications, in particular in terms of the idea of symmetry, is yet to be examined.

\section{Open groupoids and quantales}\label{sec:og}

The topology of a topological open groupoid $G$ is, similarly to that of an étale groupoid, closed under pointwise multiplication of open sets, and thus it is a quantale. A similar situation exists for open localic groupoids.
The axioms that provide the algebraic description of étale groupoid quantales in \cite{egatq} can be weakened so as to provide a characterization of the quantales (no longer unital) associated to open groupoids. Such quantales are an algebraic counterpart of open groupoids that generalizes the role played by inverse semigroups.

A continuous representation of a topological open groupoid $G$ consists of an \emph{action} of $G$ on a bundle $p:X\to G_0$ (with open $p$) \ie, a map
\[\alpha: X\times_{G_0} G_1\to X\]
satisfying suitable conditions, where $X\times_{G_0} G_1$ is the pullback of $p$ and the domain map $d$. It can be verified that $\alpha$ is necessarily an open map. This fact leads to 
an action of the open subsets of $G_1$ on those of $X$ that makes $\topology(X)$ a $\topology(G_1)$-module, and an analogous situation exists for localic groupoids.

Not surprisingly, it follows that the continuous representations of an open localic groupoid $G$ can be identified with certain modules over $\opens(G)$. Perhaps more surprisingly, the morphisms of groupoid representations can be identified with module homomorphisms (this is not true for groupoids themselves and their quantales, as we have remarked in \S\ref{sec:qg}), yielding a dual equivalence of categories.

\section{Applications}

There are many applications of groupoids in analysis, topology, geometry, and naturally also in algebra and category theory. For instance, Lie groupoids play an important role in differential geometry, and the interplay between such groupoids and operator algebras is a large part of what is meant by noncommutative geometry \cite{Connes,Paterson,Renault}, where in general one constructs C*-algebras from locally compact groupoids that are equipped with Haar measures (such groupoids are necessarily open). Some instances of this interaction are particularly well behaved. For instance, any AF-algebra is a groupoid C*-algebra of an étale groupoid (an AF-groupoid), and the relation between the two is mediated by an inverse semigroup \cite[Ch.\ III.1]{Renault}. In another direction, in topos theory open groupoids are important due to the fundamental theorem of Joyal and Tierney \cite{JT} which states that every Grothendieck topos is equivalent to the category of continuous representations of an open localic groupoid. An immediate question is how useful a reformulation of this theorem in terms of quantale modules may be.

The original motivations behind the name ``quantale'' are also related to operator algebras \cite{then}, the idea being that the underlying ``noncommutative space'' of a noncommutative C*-algebra should be a quantale, generalizing the fact that the spectrum of a commutative C*-algebra is (the spectrum of) its locale of norm-closed ideals. In addition, it was suggested \cite{then} that such a generalization of locales could provide the context for a constructive theory of noncommutative C*-algebras. This idea has led to several notions of \emph{point} of a quantale \cite{Krum02,MuPe,PeR} in the form of suitable ``simple'' modules, and to a representation theory in terms of which the ``points of noncommutative spaces'' can be classified. For instance, the equivalence classes of irreducible representations of a unital C*-algebra $A$ can be identified (albeit nonconstructively) with the points of the quantale $\Max A$ of norm-closed linear subspaces of $A$ \cite{MuPe}, a somewhat surprising consequence of this being that the quantale valued functor $\Max$ is a complete invariant of unital C*-algebras \cite{KR}. Another example is the quantale of Penrose tilings of the plane \cite{MuRe}, whose points (of a certain type) can be identified with the Penrose tilings.

Despite the progress achieved in this area, the interaction between quantales and C*-algebras is still not well understood, and attention should be given to the relations between groupoid C*-algebras, groupoid quantales, and quantales like $\Max A$.


\begin{thebibliography}{22}

\bibitem{Connes}
A.\ Connes,
Noncommutative Geometry,
Academic Press, 1994.

\bibitem{StoneSpaces}
P.T.\ Johnstone,
Stone Spaces,
Cambridge Univ.\ Press, 1982.

\bibitem{J}
P.T.\ Johnstone,
The point of pointless topology,
Bull.\ Am.\ Math.\ Soc., New Ser.\ 8 (1983) 41--53.

\bibitem{JT}
A.\ Joyal, M.\ Tierney,
An Extension of the Galois Theory of Grothendieck,
Mem.\ Amer.\ Math.\ Soc., vol.\ 309,
American Mathematical Society, 1984.

\bibitem{Krum02}
D.\ Kruml,
Spatial quantales,
Appl.\ Categ.\ Structures\
10 (2002) 49--62.

\bibitem{KR}
D.\ Kruml, P.\ Resende,
On quantales that classify C*-algebras,
Cah.\ Topol.\ Géom.\ Différ.\ Catég.\ 45 (2004) 287-296.

\bibitem{Lawson}
M.V.\ Lawson,  Inverse Semigroups --- The Theory of Partial Symmetries, World Scientific, 1998.

\bibitem{MM}
S.\ Mac Lane, I.\ Moerdijk, Sheaves in Geometry and Logic --- A First Introduction to Topos Theory, Springer-Verlag, 1992.

\bibitem{then}
C.J.\ Mulvey,
\&,
Rend.\ Circ.\ Mat.\ Palermo (2) Suppl.\
(1986) 99--104.

\bibitem{Mul:enc}
C.J.\ Mulvey,
Quantales,
in: M.\ Hazewinkel (Ed.), The Encyclopaedia of Mathematics, third supplement, Kluwer Academic Publishers, 2002, pp.\ 312--314.

\bibitem{MuPe}
C.J.\ Mulvey, J.W.\ Pelletier,
On the quantisation of points,
J.\ Pure Appl.\ Algebra
159 (2001) 231--295.

\bibitem{MuRe}
C.J.\ Mulvey, P.\ Resende,
A noncommutative theory of Penrose tilings,
Int.\ J.\ Theoret.\ Phys.\ 44 (2005) 709--743.

\bibitem{Paterson}
A.L.T.\ Paterson,
Groupoids, Inverse Semigroups, and Their Operator Algebras,
Birkh\"{a}user, 1999.

\bibitem{PeR}
J.W.\ Pelletier, J.\ Rosick\'{y},
Simple involutive quantales,
J.\ Algebra 195 (1997) 367--386.

\bibitem{Renault}
J.\ Renault, A Groupoid Approach to C*-algebras, Lect. Notes Math. 793, Springer-Verlag, 1980.

\bibitem{egatq}
P.\ Resende,
\'{E}tale groupoids and their quantales,  Preprint, 2004; arXiv:math/0412478.

\bibitem{Re04}
P.\ Resende,
Sup-lattice 2-forms and quantales,  J.\ Algebra 276 (2004) 143--167.

\bibitem{Rosenthal}
K.\ Rosenthal,
Quantales and Their Applications,  Pitman Research Notes in Mathematics Series 234, Longman Scientific \& Technical, 1990.

\bibitem{Weinstein}
A.\ Weinstein,
Groupoids: unifying internal and external symmetry, Notices Amer.\ Math.\ Soc.\ 43 (1996) 744--752.

\end{thebibliography}
\end{document}